\documentclass[review]{elsarticle}

\usepackage{lineno,hyperref}
\usepackage[tbtags]{amsmath}
\usepackage{float}
\usepackage{amssymb}
\usepackage{amsthm}
\usepackage{nomencl}
\usepackage{color}
\usepackage{amsmath}
\makeatletter
\def\BState{\State\hskip-\ALG@thistlm}
\makeatother
\usepackage[noend]{algpseudocode}
\usepackage{csquotes}
\usepackage{color}
\usepackage{rotating}
\usepackage{graphicx}
\usepackage{tikz}
\usetikzlibrary{shapes,arrows}
\usepackage{pgfplots}
\usepackage{textcomp}
\usepackage{siunitx}
\usepackage[english,ruled,vlined]{algorithm2e}
\pgfplotsset{width=10cm,compat=1.9}
\graphicspath{ {Figures/} }
\usepackage{subfigure,booktabs}
\usepackage{etoolbox}
\usepackage{multirow}

\renewcommand\nomgroup[1]{%
    \item[\bfseries
    \ifstrequal{#1}{A}{Sets}{%
    \ifstrequal{#1}{B}{Parameters}{%
    \ifstrequal{#1}{C}{Binary Variables}{%
    \ifstrequal{#1}{D}{Continuous Variables}{}}}}%
]\vspace{0.1in}}
\usepackage[top = 1in, bottom = 1in, left = 0.7in, right = 0.7in]{geometry}

\newcommand{\etal}{\textit{et al}. }
\newcommand{\ie}{\textit{i}.\textit{e}., }

\allowdisplaybreaks

\modulolinenumbers[5]

\journal{Computers \& Industrial Engineering}

\begin{document}

\begin{frontmatter}

\title{On Scheduling a Photolithography Process Containing Cluster Tools}


\author[mymainaddress]{Sreenath Chalil Madathil\corref{mycorrespondingauthor}}
\cortext[mycorrespondingauthor]{Corresponding author}
\ead{schalil@g.clemson.edu}

\author[mysecondaryaddress]{Siddhartha Nambiar}
\author[mymainaddress]{Scott J. Mason}
\author[mymainaddress]{Mary E. Kurz}

\address[mymainaddress]{110 Freeman Hall , Clemson University , Clemson, SC}
\address[mysecondaryaddress]{400 Daniels Hall , North Carolina State University , Raleigh, NC}

\begin{abstract}
Photolithography is typically the bottleneck process in semiconductor manufacturing. In this paper, we present a model for optimizing the scheduling of the photolithography process in the presence of both individual and cluster tools. The combination of these individual and cluster tools that process various layers (stages) of the semiconductor manufacturing process flow is a special type of flexible flowshop. We seek separately to minimize total weighted completion time and maximize on-time delivery performance. Experimental results suggest that our solution algorithms show promise for real world implementation as they can help to improve resource utilization, reduce job completion times, and decrease unnecessary delays in a wafer fab. 
\end{abstract}


\end{frontmatter}


\section{Introduction}
Scheduling and sequencing are indispensable processes in industry. A well-designed scheduling system helps the industry focus on increasing throughput by reducing the run time of machines, thereby saving money. Processing jobs on a ``first-come, first-serve'' basis may not be an optimal policy on the factory floor \cite{conway2012theory}. The semiconductor wafer fabrication industry is one of the largest industrial manufacturing segments. Implementing a proper scheduling system in wafer fabrication can help increase profit margins as well as reduce the time required to produce the wafers that contain integrated circuits. 

In semiconductor manufacturing, photolithography is normally one of the bottleneck processes that require high capital investments \cite{sha2006dispatching}. Hence, optimizing the photolithography process by efficiently scheduling the jobs could be beneficial for the industry. Machines that perform various steps in photolithography can be organized as a flexible flowshop system. A flexible flowshop is defined as a system in which the jobs need to be processed at different sequential stages and at least one of the stages has more than one machine operating in parallel. With the advancement of technology and because of their efficiency and profitability, cluster tools were added to the wafer fabrication processes in recent years. A cluster tool combines various types of machines that perform individual processes and organizes them around a robotic wafer transport device \cite{yim1999scheduling}. These tools consist of those machines that are capable of processing two or more stages and combine several processing modules into a single machine \cite{lee2008review}.

In this research, we develop a scheduling model for the photolithographic process, which is a special type of flexible flowshop (FFS) that has cluster tools along with the traditional individual photolithography tools. {\color{blue} According to Chiang \cite{chiang2013enhancing}, photolithography scheduling is more complex than tradition flexible flowshop scheduling. The author reviews several  reasons for this scheduling complexity such as re-entrant job flow, a jobs' readiness, due dates,  multiple machine types, multiple orders per job, and lot priorities.} Each of the jobs that enter the system typically re-visits equipment visited at earlier manufacturing (\ie reentrant flow). If the proposed model is tested successfully, it could be implemented in the semiconductor industry that employs photolithography machines with advanced cluster tools. Wafer fabs will be able to schedule their machines to improve utilization of the machines, reduce the processing time for jobs, and efficiently schedule without introducing unnecessary delays in the process.

\section{Literature Review}

Most manufacturing industries face various challenges such as processing high priority jobs, unforeseen breakdowns, scheduled maintenance, delayed processing of jobs, and meeting deadlines set by customers. Proper production planning and process scheduling help to maintain or improve the efficiency of systems and control of operations \cite{pinedoscheduling}. The significance of proper production scheduling comes to light in this scenario when manufacturers need to satisfy customer demands with the help of a minimal number of photolithography tools missing no committed completion time. This committed completion time is the due date \cite{pinedoscheduling}.  Montazeri et al. explained and reviewed different scheduling rules, such as static and dynamic rules \cite{montazeri1990analysis}. Static and dynamic rules depend on the time when the rule is applied. Static rules, applied at the start of the scheduling period, have a fixed schedule and dynamic rules change as the time progress. The authors also reviews various scheduling rules, compares their performance measures for different environments and conclude that performance evaluation depends on the objective under consideration \cite{montazeri1990analysis}. 

The four basic processes involved in manufacture of integrated circuits are wafer fabrication, wafer probe, assembly and packaging, and final testing \cite{uzsoy1992review}. A wafer fabrication process includes complex procedures and technologies that involve high capital investments. The proper utilization of wafer fabs can lead to increased profit for a semiconductor wafer fabricator. Each time a wafer passes through photolithography, a new layer of required circuitry is formed on the wafer. For most wafers there will be at least 25 such layers. Since the photolithography process is repeated during wafer fabrication, overall performance of the systems is improved by improving the photolithography process \cite{arisha2004intelligent}. The high capital cost of the photolithography tools forces the wafer manufacturers to streamline the processes to utilize these machines to the fullest possible extent.
 
There are many literatures and textbooks that explains the machine environments like a single machine, parallel machines, flowshops, job shops, flexible flowshops, and flexible job shops found in industries \cite{pinedoscheduling}. Many mixed-integer programming (MIP) models for scheduling FFS are explained in \cite{sawik2011scheduling}. The book considers various scenarios of flowshop modeling with multiple machines in each stage and finite or infinite buffers between each stage. {\color{blue}According to Floudas and Lin \cite{floudas2005mixed}, many scheduling problems use Mixed Integer Linear Programming (MILP) to find solutions due to their rigorousness, resilience, and flexible design capabilities.} Indeed, the use of MIP models is rather popular in this regard. 

Ruiz discusses the various solution approaches for the FFS problems, which includes exact methods, heuristics, and meta-heuristics \cite{ruiz2010hybrid}. In exact methods approaches such as branch and bound, algorithms solve problems to optimality. The problem with branch-and-bound algorithms is that they utilize a large amount of computer processing resources and are able to solve only problems with a few jobs and stages. Often, they are also deemed to be too complex for real world problems. Lowe et al. \cite{lowe2016integrated} proposed a deterministic MIP model to schedule weekly production quantities for semiconductor manufacturing in order to meet forecasted demand over a six-month planning horizon. {\color{blue} MIP models are proposed in \cite{sawik2012batch} for deterministic batch or cyclic scheduling in flow shops with parallel machines and finite in-process buffers. Further, \cite{sawik2014mixed} presented a new MIP formulation for cyclic scheduling in flow lines with parallel machines and finite in-process buffers, where a Minimal Part Set (MPS) in the same proportion as the overall production target is repetitively scheduled.}

A simple, two-stage flexible flowshop is strongly NP-hard \cite{hoogeveen1996preemptive}. {\color{blue} According to Kyparisis and Koulamas \cite{kyparisis2001note}, minimizing total weighted completion time for a multiple stage flexible flowshop scheduling problem is NP hard.} Hence by extension, the complexity of scheduling a larger flexible flowshop with multiple machines in almost every stage of its processing is also strongly NP hard. When compared to traditional flowshops, a photolithography system involving cluster tools, constraints for multiple wafer routes, reentrant flow, and no buffers inside the cluster tool are therefore also strongly NP hard \cite{yim1999scheduling}. Since the practical-sized complex FFS problems NP-hard, we require smart heuristics to arrive at good solutions \cite{jungwattanakit2007heuristic}.

Solving FFS problems by heuristic methods like dispatching rules and variants of shifting the bottleneck procedure (SBP) \cite{cheng2001shifting} are explained by \cite{lee2008review}. Sarin et al. \cite{sarin2011survey} provides an overview of advanced dispatching rules and compares the effectiveness of the performance from various simulation studies in a wafer fab. These dispatching rules include scheduling of general wafer fab, specific operations at bay level like photolithography, batch processing etc. The primary characteristics that make wafer fab scheduling such a different problem includes batching, reentrant flow, sequence dependent setups,  and parallel machines \cite{monch2011survey}. 

Dispatching rules include certain rules of thumb for the priority assignment of jobs onto machines. Some examples of dispatching rules include Shortest Processing Time (SPT), Longest Processing Time (LPT), and Shortest Remaining Processing Time (SRPT). The SBP uses a divide-and-conquer strategy and has been proven very effective when used in combination with exact methods for solving problems. The scheduling of a flexible flowshop with cluster tools is performed via simulated annealing \cite{yim1999scheduling} to obtain a near-optimal solution. However, the study does not consider the re-entry of jobs to previous stages. {\color{blue}Pan et al. provide a recent comprehensive literature review of the scheduling of cluster tools in semiconductor manufacturing \cite{pan2017scheduling}.}

Genetic Algorithms (GA) are a popular tool used in a number of papers focused on applications in real-world problems \cite{oduguwa2005evolutionary}. GAs have been adapted to solve problems involving sequence-dependent setup times, several production stages with unrelated parallel machines at each stage, and machine eligibility \cite{ruiz2006geneticalgorithm}. The choice of how the GA solution is represented is an important facet in the design of a GA, as representation affects other design choices, such as crossover and mutation functions. A commonly employed representation scheme is the topological ordering of the tasks. Ramachandra and Elmaghraby \cite{ramachandra2006sequencing} minimize the weighted sum of completion times in a flexible flowshop by representing the chromosomes as topological orderings of jobs, the schedules of which are obtained using a first-available machine rule for machine assignments. 


Even though most of the papers reviewed have mentioned either the scheduling of flowshops, the scheduling of flexible flowshops, and/or scheduling of cluster tools separately, there exist no efficient models that analyze a flexible flowshop that contains cluster tools and reentrant job flow across multiple product types. We will also consider job ready times and the continuous flow of jobs inside cluster tools. In this research, we develop a scheduling model for the photolithographic process, that has cluster tools along with traditional photolithography tools, and considers reentrant job flow across multiple product types. Additionally, we use two heuristic algorithms to provide numerical results. 

\begin{table}[ht!]
\centering
\color{blue}
\caption{Summary of Relevant Literature}
\label{Tab:Ref}
\begin{tabular}{lllll}
\toprule
\textbf{Reference} & \textbf{Objective} & \textbf{Methodology} & \textbf{Features} & \textbf{Complexity}\\
\cmidrule{1-5}
\cite{yim1999scheduling} & $C_{max}$ & Heuristic & Re-entry jobs &NP hard\\
 &  & Simulated Annealing & No stand-alone tools  &\\
 &  & Candidate list & No ready times  &\\
\cite{kyparisis2001note} & $\sum W_jCj$ & Approximation & Parallel machines & NP hard \\
                         &                   & WSPT heuristic & No ready time  & \\
                         &                   & Worst-case analysis & No cluster tools  &\\    
\cite{kock2007cycle}     &  $C_{max}$     & Discrete event simulation & No ready time & N/A \\
                         &  Throughput    &  Aggregate modeling approach & no re-entry &\\
\cite{zhou2014scheduling} &  $C_{max}$    &   Heuristic           & Re-entrant flow &NP hard\\  
                            &            &                       & Cluster tools & \\
\cite{park2017models} & $C_{max}$ & Fab-level simulation & Lot cycle-time & N/A  \\
                         &                   &  & Lot residency time &  \\
                         &                   &  & Lot throughput time & \\
\cite{zhang2017improved} &$C_{max}$ & Heuristic & Ready times & NP hard \\                        & & Imperialist competitive algorithm & Re-entry jobs &\\
         & & Rolling horizon & No cluster tools &\\
\bottomrule
\end{tabular}
\end{table}
\section{Problem Description}
\begin{figure}[ht!]
\centering
\includegraphics[width=\textwidth]{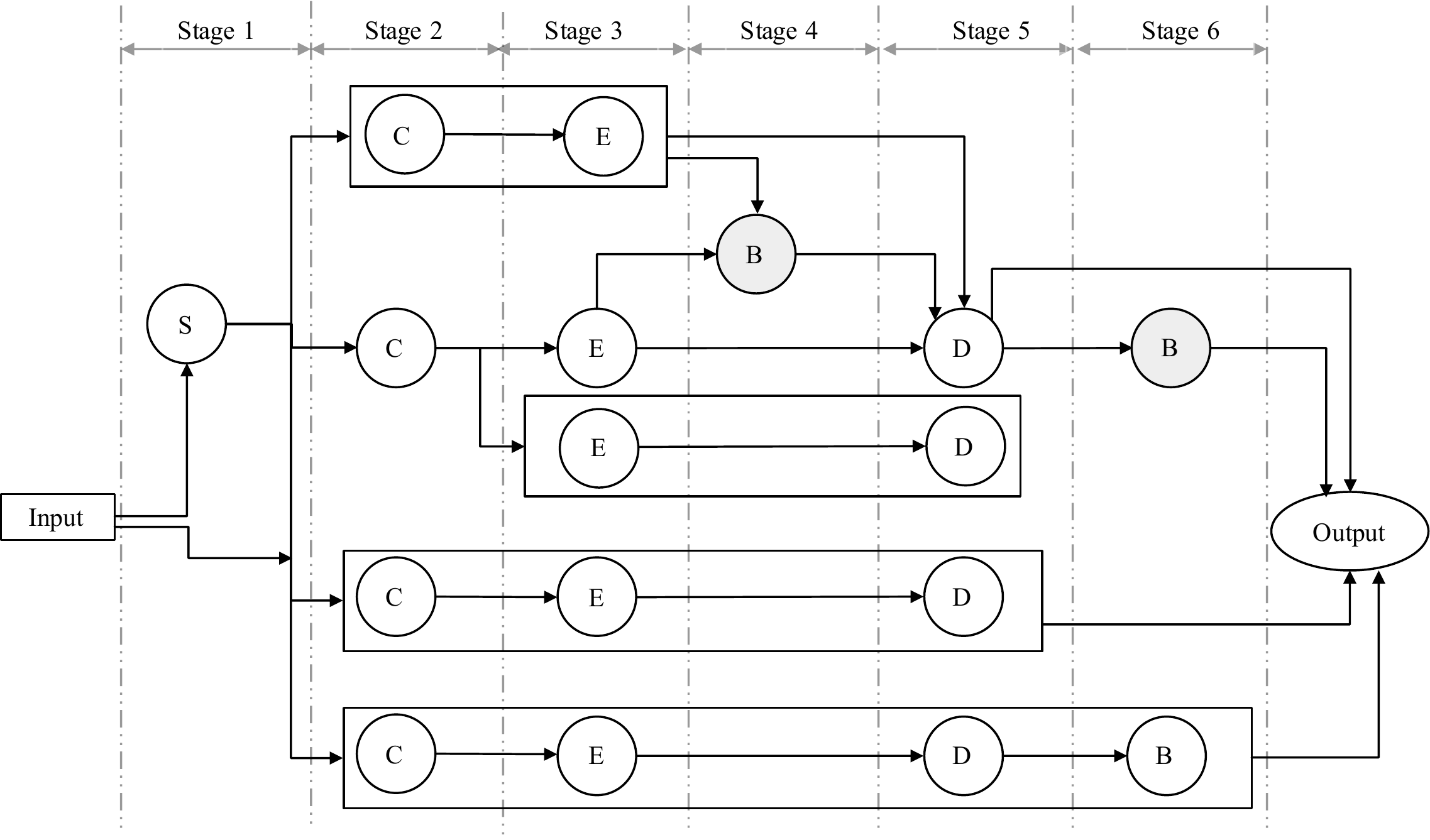}
\caption{Schematic Diagram for the Photolithography Process}
\label{Fig:photolith}
\end{figure}

The photolithography FFS system is arranged in such a way that the individual machines at each stage are organized as a general FFS with a few sets of cluster tools included. As jobs routed through the various stages of the photolithography process could return to one or more of these stages during their processing path, photolithography is a reentrant flexible flowshop \cite{graves1983scheduling}. The multistep ``photo'' process is now described in detail and is shown in Figure \ref{Fig:photolith}. The baking process in stage 4 and stage 6 can use the same bake ovens (tools) and hence represent a reentrant stage (\ie reentrant flow).

\subsection{Photolithography Stages}
In the first step of the photolithography process, a semiconductor wafer may be cleansed in a sink (tool ``S'' in Figure \ref{Fig:photolith}) \cite{mcguigan1992modeling}. The wafer is coated (tool ``C'' in Figure \ref{Fig:photolith}) with photosensitive resist and is exposed (``E'' in Figure \ref{Fig:photolith}) to light. Wafers are exposed to light with the help of a pattern mask that controls the wafer areas that receive light exposure. This helps to define the required circuit functionality. The exposed wafer is then developed (``D'' in Figure \ref{Fig:photolith}) so that the required patterns are imprinted onto the wafer by removing the exposed photoresist. The final photolithography stage is baking (``B'' in Figure \ref{Fig:photolith}). Sometimes, wafers are baked before and/or after the developing stage. 

The flow diagram of the photolithography process for a single layer of wafer fabrication is illustrated in Figure \ref{Fig:photolith}. Normally, a wafer repeats this process 20-30 times during its process flow. Figure \ref{Fig:photolith} also depicts cluster tools that are used in the photolithography process. Cluster tool ``CEDB'' processes the Coat, Expose, Develop, and Bake steps in order. Similarly cluster tools ``CED,'' ``CE'', and ``ED'' process Coat-Expose-Develop (CED), Coat-Expose (CE), and Expose-Develop (ED), respectively. A job can have multiple routes based on the process flow. For example, a job that requires a coat, expose, develop and bake processes can flow through any of the 4 paths using the combinations of tools like C$\rightarrow$E$\rightarrow$D$\rightarrow$B, CE$\rightarrow$D$\rightarrow$B, CED$\rightarrow$B or CEDB.

{\color{blue} The transportation of wafer lots inside a wafer fab (interbay and intrabay) is via automated material handling systems (AMHSs). AMHSs transport lots both to stockers and to downstream production tools, as dictated by the shop floor control system (MES), as soon as the current tool the lot is on completes its processing. The availability of the downstream tool determines the destination of completed lots, whether to the stocker or to the downstream processing tool. This wafer transportation inside a wafer fab takes place asynchronously with wafer transported to the next tool as soon as current tool completes its process and the tool for the next process is free \cite{kock2007cycle}. Another transport technique is synchronous in which transportation of the all wafers to the next process happen when the slowest wafer has finished processing. The WIP stockers in wafer fabs are sized to hold all required work-in-process within the fab \cite{cardarelli1995simulation}. In this way, lots waiting to be processed do not wait in a physical queue located on/near the production equipment, but rather in a virtual queue as house in the WIP stocker. Stockers are typically located throughout the wafer fab, often at the end of each processing bay along the main center aisle.}

\subsection{Methodology}
{\color{blue} We propose a MILP model to schedule this photolithography stages and attribute our model's combinatorial nature to various discrete decisions such as job's assignment to machines and sequencing of jobs in a machine. We develop two heuristic algorithms, compare the solution quality and solving time for the three methods, and evaluate the heuristic's performance for larger problems.} 

The set of jobs entering the system for processing can be characterized by their ready times ($r_j$), the time at which the job is released to the shop by some external job scheduler \cite{conway2012theory}. An initial step in optimizing industrial processes often includes improving the total execution time of machines, also known as makespan ($C_{max}$). A schedule that minimizes makespan can be obtained by applying MIP techniques \cite{floudas2005mixed}. 

In this research, we develop a MIP formulation to optimize the makespan of the photolithography process containing reentrant jobs with ready times. In terms of the standard $\alpha | \beta | \gamma$ scheduling notation introduced in Graham \etal \cite{graham1979optimization}, the problem under consideration is defined as $FFm | r_j,rcrc| C_{max}$ ~\cite{pinedoscheduling}. Other objectives that could be optimized in a scheduling system include total weighted completion time (WCT) and total weighted tardiness (TWT). The total weighted completion time is represented as $\sum w_kC_k$ with $w_k$  denoting the weight or priority of job $k$ and $C_k$ representing the completion time of job $k$. In practice, WCT is a surrogate measure of the inventory or holding cost incurred by the schedule \cite{pinedoscheduling}. The total weighted tardiness, $\sum w_kT_k$  where $T_k$  is the tardiness of the job $k$, is generally an objective that relates to on-time delivery.


\section{Mathematical Formulation}
This section explains the model formulation along with the notations used in the model, followed by an explanation of the constraint sets. The MIP formulation for scheduling flowshops with parallel machines and infinite in-process buffers \cite{sawik2011scheduling} is used as the base model for this FFS and additional constraints are added to expand this formulation to incorporate the scheduling of cluster tools and reentrant flow.  
\subsection{Model Parameters and Variables}
\subsection {Notation}
\begin{tabular}{l l}
	\multicolumn{2}{l}{Sets}\\
	$I$&set of processing stages indexed by $i \in I = \{1,..,m\}$\\
	$J_i$&set of processors in each stage $i$ indexed by $j \in J_i = \{1,..,m_i\}$\\
	$K$&set of jobs that needs to be processed indexed by $k \in K = \{1,..,n\}$\\
	$CEDB$&set of cluster machines for $C-E-D-B$ indexed by $i_1$\\
	$CED$&set of cluster machines for $C-E-D$ indexed by $i_2$\\
	$CE$&set of cluster machines for $C-E$ indexed by $i_3$\\
	$ED$&set of cluster machines for $E-D$ indexed by $i_4$\\
\end{tabular}

\begin{tabular}{l l}
	\multicolumn{2}{l}{Parameters}\\
	$m$&number of processing stages\\
	$m_i$&number of machines at each processing stage i\\
	$n$&number of jobs\\
	$P_{ik}$&processing time for job $k$ in stage $i$\\
	$r_k$&ready time for job $k$\\
	$d_k$&due date for job $k$\\
	$w_k$&priority for job $k$\\
	$M$&sum of the processing time of all jobs in the system\\	
\end{tabular}

\begin{tabular}{l l}
	\multicolumn{2}{l}{Decision Variables}\\
	$C_{ik}$&completion time of job $k$ at stage $i$\\
	$C_{max}$&makespan or time at which all jobs complete their operations on all stages\\
	$x_{ijk}$&1, if job $k$ is assigned to machine $j$ in stage $i$\\
	&0, otherwise\\
	$y_{kl}$&1, if job $k$ precedes job $l$ in the processing sequence\\
	&0, otherwise\\	
\end{tabular}
\subsection{Model Formulation}
\vspace{-1cm}
\begin{align}
\label{eq01}&\min ~ C_{max}\\
\label{eq02}&{C_{1k}} \ge {P_{1k}} + {r_k} & & ~\forall~ k \in K\\
\label{eq03}&{C_{ik}} - {C_{\left( {i - 1} \right)k}} \ge {P_{ik}} & & ~\forall~ k \in K,i \in I,i > 1\\
\label{eq04}&{C_{ik}} + M\left( {2 + {y_{kl}} - {x_{ijk}} - {x_{ijl}}} \right) \ge {C_{il}} + {P_{ik}} & &~\forall~ i \in I,j \in {J_i},k \in K,l \in K,l > k\\
\label{eq05}&{C_{ik}} + M\left( {3 - {y_{kl}} - {x_{ijk}} - {x_{ijl}}} \right) \ge {C_{ik}} + {P_{il}} & &~\forall~ i \in I,j \in {J_i},k \in K,l \in K,l > k\\
\label{eq06}&{C_{mk}} \le {C_{\max }}  & &~\forall~ k \in K\\
\label{eq07}&\sum\limits_{j \in {J_i}} {{x_{ijk}}}  = 1 & &~\forall~ i \in I,k \in K,P_{ik} > 0\\
\label{eq08}&\sum\limits_{j \in {J_i}} {{x_{ijk}}}  = 0 & &~\forall~ i \in I,k \in K,{P_{ik}} = 0\\
\label{eq09}&\sum\limits_{i \in I} {{x_{i{i_1}k}}}  = 3{x_{2{i_1}k}} & &~\forall~ {i_1} \in {\rm{CEDB}},k \in K,{i_1} \ne 0, 3 \le i \le m,i \ne 4\\
\label{eq10}&{x_{4jk}} \le 1 - {x_{2{i_1}k}}  & &~\forall~ {i_1} \in {\rm{CEDB}} \cup {\rm{CED}},{i_1} \ne 0,j \in {J_4},k \in K\\
\label{eq11}&{C_{2k}} + M(2 + {y_{kl}} - {x_{2{i_1}k}} - {x_{2{i_1}l}}) \ge {C_{ml}} + {P_{2k}} & &~\forall~ {i_1} \in {\rm{CEDB}},{i_1} \ne 0,k \in K,l \in K,l > k\\
\label{eq12}&{C_{2l}} + M(3 - {y_{kl}} - {x_{2{i_1}k}} - {x_{2{i_1}l}}) \ge {C_{mk}} + {P_{2l}} & &~\forall~ {i_1} \in {\rm{CEDB}},{i_1} \ne 0,k \in K,l \in K,l > k\\
\label{eq13}&\sum\limits_{i \in I} {{x_{i{i_2}k}}}  = 2{x_{2{i_2}k}} & & ~\forall~ {i_2} \in {\rm{CED}},k \in K,{i_2} \ne 0,3 \le i \le m - 1,i \ne 4\\
\label{eq14}&{C_{2k}} + M(2 + {y_{kl}} - {x_{2{i_2}k}} - {x_{2{i_2}l}}) \ge {C_{(m - 1)l}} + {P_{2k}}& &  ~\forall~ {i_2} \in {\rm{CED}},{i_2} \ne 0,k \in K,l \in K,l > k\\
\label{eq15}&{C_{2l}} + M(3 - {y_{kl}} - {x_{2{i_2}k}} - {x_{2{i_2}l}}) \ge {C_{(m - 1)k}} + {P_{2l}}& &  ~\forall~ {i_2} \in {\rm{CED}},{i_2} \ne 0,k \in K,l \in K,l > k\\
\label{eq16}&\sum\limits_{i \in I} {{x_{3{i_3}k}}}  = {x_{2{i_3}k}}& & ~\forall~ {i_3} \in {\rm{CE}},k \in K,{i_3} \ne 0\\
\label{eq17}&{C_{2k}} + M(2 + {y_{kl}} - {x_{2{i_3}k}} - {x_{2{i_3}l}}) \ge {C_{3l}} + {P_{2k}}& & ~\forall~ {i_3} \in {\rm{CE}},{i_3} \ne 0,k \in K,l \in K,l > k\\
\label{eq18}&{C_{2l}} + M(3 - {y_{kl}} - {x_{2{i_3}k}} - {x_{2{i_3}l}}) \ge {C_{3k}} + {P_{2l}}& & ~\forall~ {i_3} \in {\rm{CE}},{i_3} \ne 0,k \in K,l \in K,l > k\\
\label{eq19}&\sum\limits_{i \in I} {{x_{5{i_4}k}}}  = {x_{3{i_4}k}} & &~\forall~ {i_4} \in {\rm{ED}},k \in K,{i_4} \ne 0\\
\label{eq20}&{C_{3k}} + M(2 + {y_{kl}} - {x_{3{i_4}k}} - {x_{3{i_4}l}}) \ge {C_{5l}} + {P_{3k}}& & ~\forall~ {i_4} \in {\rm{ED}},{i_4} \ne 0,k \in K,l \in K,l > k\\
\label{eq21}&{C_{3l}} + M(3 - {y_{kl}} - {x_{3{i_4}k}} - {x_{3{i_4}l}}) \ge {C_{5k}} + {P_{3l}} & & ~\forall~ {i_4} \in {\rm{ED}},{i_4} \ne 0,k \in K,l \in K,l > k\\
\label{eq22}&{C_{4k}} + M(2 + {y_{kl}} - {x_{4ik}} - {x_{6il}}) \ge {C_{6l}} + {P_{4k}}& & ~\forall~ i \in {{\rm{J}}_{\rm{4}}} \cap {J_6},i \ne 0,k \in K,l \in K,l > k\\
\label{eq23}&{C_{6l}} + M(3 - {y_{kl}} - {x_{4ik}} - {x_{6il}}) \ge {C_{4k}} + {P_{6l}} & & ~\forall~ i \in {{\rm{J}}_{\rm{4}}} \cap {J_6},i \ne 0,k \in K,l \in K,l > k\\
\label{eq24}&{C_{6k}} + M(2 + {y_{kl}} - {x_{6ik}} - {x_{4il}}) \ge {C_{4l}} + {P_{6k}} & &~\forall~ i \in {{\rm{J}}_{\rm{4}}} \cap {J_6},i \ne 0,k \in K,l \in K,l > k\\
\label{eq25}&{C_{4l}} + M(3 - {y_{kl}} - {x_{6ik}} - {x_{4il}}) \ge {C_{6k}} + {P_{4l}}& & ~\forall~ i \in {{\rm{J}}_{\rm{4}}} \cap {J_6},i \ne 0,k \in K,l \in K,l > k\\
\label{eq26}&{C_{\max }} \ge 0\\
\label{eq27}&{C_{ik}} \ge 0& & ~\forall~ i \in I,k \in K\\
\label{eq28}&{x_{ijk}} \in \left\{ {0,1} \right\} & &~\forall~ i \in I,j \in {J_i},k \in K\\
\label{eq29}&{y_{kl}} \in \left\{ {0,1} \right\} & &~\forall~ k \in K,l \in K
\end{align}
The model's objective function (\ref{eq01}) minimizes makespan. The objective could be changed to minimizing the weighted completion time (WCT) if desired. Constraint sets (\ref{eq02}) and (\ref{eq03}) ensure that a job starts processing at stage 1 and processes successively on all downstream machines. The overlapping of more than one job on a single machine at a time is prevented by constraint sets (\ref{eq04}) and (\ref{eq05}). These constraints act as ``either-or'' constraints, which imply that one of the constraints will be active for a particular value of $y_{kl}$ . When job $k$ precedes job l, then constraint set (\ref{eq04}) will be active and constraint set (\ref{eq05}) will be inactive, because of the value of ``M'' and vice-versa. The value of M is assigned as the sum of processing time of all jobs in the system. The maximum completion time, $C_{max}$, should be greater than or equal to the completion time of the last job in the final stage of the processing. This is achieved by including the constraint sets (6). Constraint sets (\ref{eq07}) and (\ref{eq08}) ensure that if the processing time for a job in any stage is a non-zero, positive number, then one of the machines in that stage must process the job \cite{fourer1990modeling}. If the processing time for a job at any stage is zero, then that stage is skipped for that particular job.

The remaining sets of constraints are developed for cluster tools and reentrant flow processes. Constraint sets (\ref{eq09}) - (\ref{eq12}) model the cluster process of coat, expose, develop, and bake in a single machine. Constraints set (\ref{eq10}) ensures that a job that needs to be processed on stage 4 bake process will not enter the cluster tools CEDB or CED. Constraint sets (\ref{eq09}), (\ref{eq13}), (\ref{eq16}) and (\ref{eq19}), ensure that a job that enters a cluster machine will stay inside that machine until it completes all the processes performed by the cluster tool. Pairs of constraint sets (\ref{eq11}) - (\ref{eq12}), (\ref{eq14}) - (\ref{eq15}), (\ref{eq17}) - (\ref{eq18}), and (\ref{eq20}) - (\ref{eq21}) stop jobs from entering the cluster tool if the machine is already processing some other job. Constraint set (\ref{eq10}) updates the assignment variable for job k, for the first bake stage to 0, if job $k$ does not require the baking. Constraint sets (\ref{eq22}) - (\ref{eq25}) model the re-entry of jobs in the bake process of photolithography. The photolithography process under consideration has re-entry at stages 4 and 6. Hence, the machines of the fourth and sixth stages $J_4$ and $J_6$ are referenced in these equations. These constraints guarantee that if any of the jobs is being processed in the bake oven at any of its two stages, \ie the first bake process or the reentrant second baking stage, then no other job will enter the machine. Finally, constraint sets (\ref{eq26}) - (\ref{eq29}) are non-negativity constraints, which imply that these variables should have a value greater than or equal to zero.

For calculating the minimum TWT, the objective that will be used is as below.
\begin{equation} 
\begin{aligned}
\label{eq30}\text{minimize}\sum{w_kT_k}\end{aligned} 
\end{equation}
Constraint set (\ref{eq06}) will be replaced by a new constraint to incorporate the tardiness value. Constraint set (\ref{eq31}) is used when the objective function is TWT.
\begin{equation} 
\begin{aligned}
\label{eq31}C_{mk}-d_k \leq T_k ~\forall~ k \in K
\end{aligned} 
\end{equation}

\section{Experimental Study}
\subsection{Experimental Plan}
Our experimental plan evaluates three objective functions with the proposed MIP for photolithography scheduling: minimizing makespan ($C_{max}$), minimizing total weighted completion time (WCT), and minimizing total weighted tardiness (TWT). The experimental design used in the random problem generation \cite{mehta1998minimizing} is given in Table \ref{Tab:T1}. These parameter values were vetted with the industrial partner to be an appropriate set of factors for experimentation. Three different levels of the number of jobs to be scheduled are investigated in this problem:  5, 10, and 25. Each set of jobs tested with two levels of ready times. For the first condition, all jobs have zero ready time and for the second condition, some portion of the jobs have a ready time that is a non-zero, randomly generated value while the remaining jobs of the same sets have zero ready time.
\begin{table}[ht!]
\centering
\caption{Experimental Design}
\label{Tab:T1}
\begin{tabular}{llll}
\toprule
\textbf{Experimental Factor} & \multicolumn{3}{l}{\hspace{2.5cm} \textbf{Settings}} \\
\cmidrule{1-4}
Number of jobs, $n$ &\multicolumn{3}{l}{ 5     \hspace{2.5 cm} 15            \hspace{2.5 cm}    20} \\
Ready time & $r_k$ = 0 for all $k$ & \multicolumn{2}{l}{30\% of jobs,  $r_k$ = 0} \\
& & \multicolumn{2}{l}{70\% of jobs,  $r_k$ = RANDOM{[}1,$^2/_3 \times C_{max}${]}} \\
Job due date $d_k$ &  \multicolumn{3}{l}{$d_k$ = Uniform {[}$\mu$(1-0.5$R$) , $\mu$(1+0.5$R$){]}} \\
& \multicolumn{3}{l}{with} \\
& \multicolumn{3}{l}{$\mu $ = $C_{max} \times (1-T)$} \\
& \multicolumn{3}{l}{$C_{max} $ = 1.5 $\times$ ($n \times \frac{P_{BN}} {m_{iBN}}+ P_{NBN}$)} \\
& \multicolumn{3}{l}{$P_{BN} $ = Processing time of the bottleneck stage} \\
& \multicolumn{3}{l}{$m_{iBN} $ = Number of machines that processes the bottleneck stage} \\
& \multicolumn{3}{l}{$P_{NBN} $ = Sum of the Processing time of all other non-bottleneck stages} \\
& \multicolumn{3}{l}{$T $ = 0.3 and 0.6} \\
& \multicolumn{3}{l}{$R $ = 0.5 and 2.5} \\
Processing time & Stage 1 & \multicolumn{2}{l}{80\% of jobs $P_{1j}$ = 40} \\
 &  & \multicolumn{2}{l}{20\% of jobs $P_{1j}$ = 0} \\
 & Stage 2 & \multicolumn{2}{l}{100\% of jobs $P_{2j}$ = 20} \\
 & Stage 3 & \multicolumn{2}{l}{100\% of jobs $P_{3j}$ = 75} \\
 & Stage 4 & \multicolumn{2}{l}{20\% of jobs $P_{4j}$ = 45} \\
 &  & \multicolumn{2}{l}{80\% of jobs $P_{4j}$ = 0} \\
 & Stage 5 & \multicolumn{2}{l}{100\% of jobs $P_{5j}$ = 30} \\
 & Stage 6 & \multicolumn{2}{l}{50\% of jobs $P_{6j}$ = 45} \\
 &  & \multicolumn{2}{l}{50\% of jobs $P_{6j}$ = 0} \\
\bottomrule
\end{tabular}
\label{Tab:ExpDes}
\end{table}

The due date value is generated using a discrete uniform distribution (Table \ref{Tab:T1}). The calculation of the estimated makespan includes the total number of jobs, the processing time of the photolithography process's bottleneck stage, the total number of machines that process the bottleneck stage, and the sum of the processing times for all non-bottleneck stages in photolithography. The parameter T is the expected percentage of tardy jobs and two cases for the value of T are considered in this experimentation:  0.3 and 0.6. Further, R is a range parameter that is studied at two levels:  0.5 and 2.5 \cite{mehta1998minimizing}.

The weights ($w_k$) are calculated based on a random distribution of all integers between 1 and 5, 1 being a low priority job and 5 being a high priority job. Considering the three levels for the number of jobs, two scenarios for job ready times, and four combinations of due date parameters T and R, 24 unique combinations of data are run for two different resource (equipment) levels (Table \ref{Tab:T2}). The equipment settings in scenario 1 are based on representative data obtained from an actual wafer fab that partnered with the authors to conduct this study. In order to examine potential performance differences, scenario 2 was created by reducing equipment counts from scenario 1. As three different objective functions are investigated, a total of 24(2)(3) = 144 unique scenarios exist for investigation. Based on 10 replications for each unique scenario, a total of 1,440 files are generated for analysis and comparison. Microsoft Excel 2010 is used to generate random numbers for the various cases in the experimental plan.
\begin{table}[ht!]
\centering
\caption{Equipment settings}
\label{Tab:T2}
\begin{tabular}{lll}
\toprule
\textbf{Machine Type} & \textbf{Scenario 1} & \textbf{Scenario 2} \\
\cmidrule{1-3}
Sink                  & 4                   & 2                   \\
Coat                  & 2                   & 1                   \\
Expose                & 4                   & 2                   \\
Develop               & 2                   & 1                   \\
Bake                  & 3                   & 2                   \\
CE                    & 2                   & 1                   \\
CED                   & 2                   & 1                   \\
CEDB                  & 2                   & 1                   \\
ED                    & 1                   & 1                 \\
\bottomrule
\end{tabular}
\end{table}

\subsection{Model Execution}
The random data that was created is used to test the mathematical model. After implementing the model in AMPL \cite{fourer1993ampl}, 1338 test files were run and the objective functions are validated for proper machine assignments. The solution was produced using Gurobi 5.1.0 solver \cite{optimization2012inc} within a 7200 CPU seconds time limit on a Windows 7 platform with Intel \textsuperscript{\textregistered} Core \textsuperscript{TM} 2 Quad CPU Q6600 @2.40 GHz with 16GB of RAM. Although Gurobi 5.1.0 did not converge to the optimal solution within the allowed 7200 seconds for most of the instances, Gurobi 5.1.0 was able to obtain a good solution quickly. {\color{blue}We explain the model size using the number of constraints and variables for each of the problems. Table \ref{Tab:NT1} provides the number of constraints and variables introduced into the model based on the number of jobs.}
\begin{table}[ht!]
\centering
\color{blue}
\caption{Model Size}
\label{Tab:NT1}
\begin{tabular}{p{0.6in}p{0.7in}p{0.7in}p{0.7in}p{0.7in}}
\toprule
\textbf{Problem Size} & \textbf{Constraints} & \textbf{Continuous Variables} & \textbf{Binary Variables}  & \textbf{Total Variables}\\
\cmidrule{1-5}
5   & 1185    & 36    &  215      & 251 \\
15  & 11205   & 106   &  795      & 901 \\
25  & 31425   & 176   &  1575     & 1751 \\
\bottomrule
\end{tabular}
\end{table}

Since the scheduling problem under study is strongly NP-Hard, heuristic and/or metaheuristic approaches may provide good, near optimal solutions that are better than the solution obtained by a time-limited MIP \cite{urlings2010heuristics}. The heuristic approach for the problem under study is presented, and then its performance is compared with that of the proposed MIP model under a time limit restriction.

\makeatletter
\def\BState{\State\hskip-\ALG@thistlm}
\makeatother
\section{Heuristic Algorithm}
A review of the available literature confirmed that no heuristic is currently available for analyzing the flexible flowshop scheduling problem with cluster tools and job ready times. Therefore, in order to very quickly obtain good solutions to the research problem under study, we created our own GA-based heuristic \cite{geneticAlgoIntro} and a constructive heuristic.

\subsection{Constructive Heuristics}
To this end, we generated a constructive heuristic (pseudocode \ref{Alg:Fit}) that takes in a topological ordering (permutation) of jobs as input and returns either $C_{max}$, TWT, or WCT using the pseudocode {\ref{Alg:CS}}. The pseudocodes \ref{Alg:Fit} and \ref{Alg:CS} illustrates this constructive heuristic.

\begin{algorithm}[ht!]
\caption{Procedure SchedulePhoto (SP)}
\SetAlgoLined
 \textbf{Initialization}\;
 Let $[O]$ denote the $o^{th}$ ordered job in a list and $n = |O|$, the cardinality of the set O\;
 Let $S_{itrcnt}$ denote the current schedule and $S_b$ denote the best schedule found so far\;
 Let $G(S_{itrcnt})$ denote the corresponding values of the objective function\;
 Let $P_{i,jcnt}$ denote the processing time of job $jcnt$ at stage $I$ and $maxitrcnt$ is the maximum number of iterations required\;
\textbf{Procedure SchedulePhoto (SP)}\;
Sort jobs in ascending order of ready times ($r_k$). Break any ties by sorting those jobs in ascending order of $\frac{d_k}{w_k}$. Finally, break any remaining ties by sorting the jobs in the descending order of total count of cluster tools that the job could be processed on.\;
Define $X = \{ k \in K | r_k = 0 \}$\;
Define $Y = \{ k \in K | 0 < r_k < \text{average ready time} \}$\;
Define $Z = \{ k \in K | r_k \geq \text{average ready time} \}$\;
\For{ $itrcnt = 1 ... maxitrcnt$} {
	\If{$itrcnt = 1$}{
    	Call \textbf{Procedure CreateSchedule (CS)} using $[O]$ and save the schedule in $S_{itrcnt}$\;
        $S_b = S_{itrcnt}$\;
    }
    \Else{
    	\If{$itrcnt \leq maxitrcnt/2$}{
        	Use a hill-climbing approach in which we randomly select two elements and swap them\;
            Update the order in $[O]$ based on swapping of elements\; 
        }
        \Else{
        	Generate a random number $U[0,1]$ for each element in $X,Y$ and $Z$ and sort each set in the ascending order of the random number\;
            Update the order in $[O]$ based on order in $X,Y$ and $Z$\;
        }
    	Call \textbf{Procedure CreateSchedule (CS)} using $[O]$ and save the schedule in $S_{itrcnt}$.
    }
    If $G(S_{itrcnt}) < G(S_b)$, then $S_b = S_{itrcnt}$\;
}
 \KwResult{Write the output results}
\label{Alg:Fit}
\end{algorithm}
\color{black}

\begin{algorithm}[ht!]
\caption{Procedure CreateSchedule (CS)}
\SetAlgoLined
\textbf{Procedure CreateSchedule (CS)}\;
\For{$ jcnt = [1] ... [n]  $}{
	\For{$i = 1 ... 6$}{
    	\If{$P_{i,jcnt} > 0$}{
        	Define set $J_{iavbl} = $ List of available machines at stage $i$ that process job $jcnt$\;
            Select a machine $m$, with the largest number of cluster tools in it, from $J_{iavbl}$\;
            Update $i$ with the number of stages processed by machine $m$\;
            Update completion time $C_{i,jcnt}$ as ($P_{i,jcnt}$ + $\max$ [ready time, completion time of previous job on that machine])\;
        }
    }
}
\KwResult{Update $S_{itrcnt}$ with the current schedule}
\label{Alg:CS}
\end{algorithm}
\color{black}

\subsection{Genetic Algorithm}

Genetic algorithms use ideas borrowed from the concepts of genetics and biological evolution. The base version of the algorithm treats solutions as genomes that are iteratively combined in pairs (binary crossover) and mutated (binary mutation) to produce offsprings in every generation. This process is repeated until a specific fitness function achieves a desired threshold criteria. In our Algorithm \ref{Alg:GA}, the fitness function computes the value of $C_{max}$, TWT, or WCT, depending on what is the desired objective function to be minimized. Further, we develop approaches to perform crossover and mutation functions that enable us to efficiently solve the flexible flowshop scheduling problem. As Algorithm \ref{Alg:GA} outlines the pseudocode for GA, we also provide a description of crossover and mutation functions in the following subsections. Our algorithm was implemented in MATLAB 8.1.0.604 (R2016a) and uses the same data as input that was generated for the mathematical model. With regards to the threshold criteria, we set the algorithm to terminate after MaxGenerations = 500 generations or if the average relative change in the best fitness function value in 50 continuous generations is less than or equal to $10^{-6}$. Finally, we also set the size of the population in each generation (popsize) to 100.

\vspace{0.5cm}
\begin{algorithm}[ht!]
\caption{Genetic Algorithm-based Heuristic}
\SetAlgoLined
 Initialization\;
 $t \,\leftarrow \, 0$\;
 Define $S$ as the set containing all permutations of $(1,2,3,...,n)$\;
 Define $P(t) \subset S$ as the initial population of permutations\;
 Define $F$ as a fitness function to evaluate the fitness of each individual  $p \in P(t)$\;
 Define $arrF$ as an array to store the values of the fitness functions\;
 \While{$t < \text{Max Generations}$}{
 	Build parent population $P_p(t)$ based on selection criterion\;
    \For{$t = 1 \text{ to } \text{ popSize}$}{
      Randomly select two elements $(u, v) \in P_p(t)$\;
      Perform crossover function $C$ to obtain offspring. $u_o = C(u, v)$\;
      Apply Mutation function $M$ on offspring\;
      \If{$F(u_o) > F_{best}$}{
      	$F_{best} \leftarrow F(u_o)$\;
        $S_{best} \leftarrow u_o$\;
      }
      Add $u_o$ to new generation's population\;
    }
    $t \leftarrow t + 1$\;
    \If{Average change in last 50 elements in $arrF < \epsilon$}{
      	$\textbf{break}$\;
      }
 }
 \KwResult{Write the output results}
\label{Alg:GA}
\end{algorithm}
\subsubsection{Crossover Function}
In genetic algorithms, crossover functions modify the makeup of genomes from one generation to the next by taking more than one parent genome and producing a child genome. In our GA presented in Algorithm \ref{Alg:GA}, we perform the crossover function in the following manner.

Let us assume two parent genomes are coded as $u$ and $v$ which are a vector of integers representing permutations of $\{1,2,3,...,n\}$. Define $u^i$ as the $i^{th}$ element of $u$ and generate a random integer $r$ from the set of integers $\{1,2,3,...,n\}$. Now, crossover is performed by swapping the positions of elements $u^r$ and $v^r$ in both $u$ and $v$. For instance, let $u = \{2,1,4,3,5\}$ and $v = \{5,2,3,4,1\}$. Let random integer $r = 2$. Therefore, we swap the positions of the second elements of $u$ and $v$ respectively, i.e., elements 1 and 2 in both $u = \{2,1,4,3,5\}$ and $v = \{5,2,3,4,1\}$. The newly formed children are now $\{1,2,4,3,5\}$ and $\{5,1,3,4,2\}$. One of these children is selected at random and added to the newly formed pool of children to be mutated in the next step. 
\subsubsection{Mutation Function}
During mutation in GA, one or more gene values are altered slightly so as to ensure that genetic variation is preserved. In Algorithm \ref{Alg:GA}, each of the children formed after the crossover step is modified slightly by taking two random elements of the permutation and switching their positions. For instance, if one of the children formed from the crossover step is $\{4,2,1,3,5\}$, two elements are picked at random, (say, $3$ and $4$) and their positions are swapped, to produce a new offspring with a slightly altered genetic makeup $: \{3,2,1,4,5\}$. 
\subsubsection{Fitness Function}
The proposed heuristic solution outlined in Algorithm \ref{Alg:GA} involves computing a fitness function. Procedure CS in Algorithm 3 (pseudocode \ref{Alg:CS}) computes the fitness function depending on the objective function such as  $C_{max}$, TWT, or WCT, that we wish to minimize.

\section{Results}
\subsection{Solution time performance}
\begin{figure}[ht]
\centering
\includegraphics[width=\textwidth]{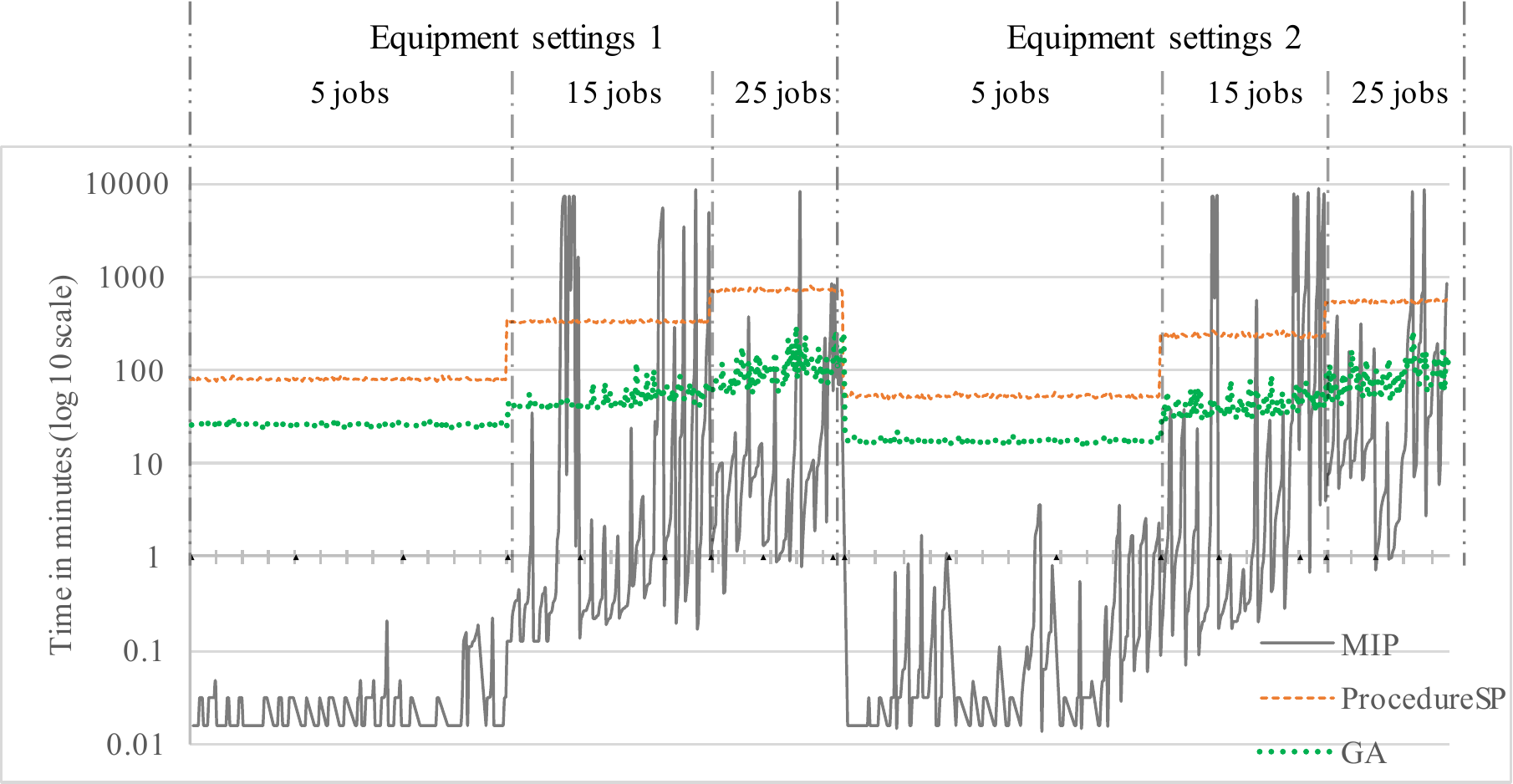}
\vspace{-1cm}\caption{Solution time comparison for optimally completed MIP instances}
\label{Fig:Optimal}
\end{figure}

\begin{figure}[ht]
\centering
\includegraphics[width=\textwidth]{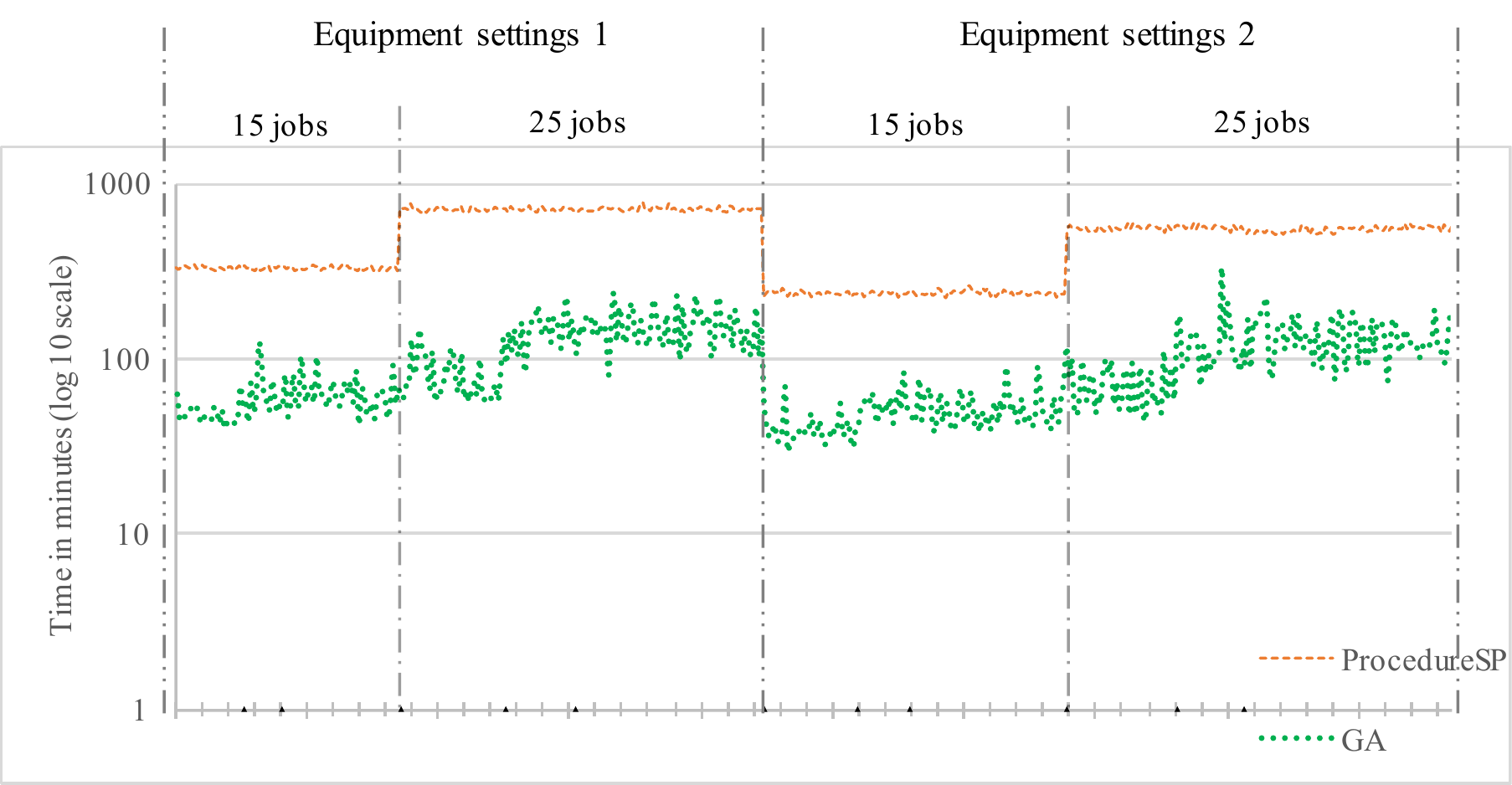}
\vspace{-1cm}\caption{Solution time comparison for time-limited MIP instances}
\label{Fig:NoOptimal}
\end{figure}

{\color{blue}We compare the algorithms' efficiency using the solving time. In order to perform a fair comparison, we classify the results in terms of the MIP's optimality status. If the MIP obtained an optimal solution, we compare the corresponding instances solution time for algorithms \ref{Alg:Fit} and \ref{Alg:GA} and is presented in Figure \ref{Fig:Optimal}. For instances with small jobs, our MIP model was working better in terms of solution time. For larger instances genetic algorithm (\ref{Alg:GA}) is better. As the size of the instances increase, the solution time also increased. The number of instances that MIP solved to optimality with 25 jobs are low when compared to the 5 job instances.}
{\color{blue}The real difference in algorithms' performances can be found when we compare the results for the timed-out MIP instances. Figure \ref{Fig:NoOptimal} compares the solution time for the two heuristics whereas the MIP was timed out after 7200 seconds. Our GA-based algorithm, consistently found better results with in certain time limit.}

\subsection{Experimental Results}

The performance of GA and the proposed mathematical model is compared by computing a performance ratio. Let performance ratio PR be defined as the ratio of the objective function value obtained by GA ($OF_{GA}$) for a problem instance to the optimal objective function value produced by the mathematical model ($OF_{opt}$) for the same problem instance (\cite{erramilli2006multiple}). While the PR ratio can be computed for any objective function case of interest, it is only valid when the MIP model produces an optimal solution. In this way, an estimate of the quality of the GA solution is obtained in terms of its percent above the optimal solution value.

Once the results are obtained for each instance, the PR values can be averaged across all experimental instances for a given type of problem type (e.g., all instances with five jobs). The set of like problem instances is characterized in terms of (n, $r_k$, T, R, mc). In this expression, n is the number of jobs and $r_k$ = 0 denotes all 0 job ready times while $r_k$ = 1 denotes the presence of non-zero ready times. Further, T and R are the due date-related parameters described above and mc represents the machine configuration (scenario 1 or scenario 2) from Table \ref{Tab:T2}. An example for this instance characterization approach is (15, 0, 0.3, 2.5, 1), which represents the average PR values for problems with 15 jobs that have zero job ready times, T and R values of 0.3 and 2.5, respectively, and machine configuration of scenario 1 from Table \ref{Tab:T2}.

\begin{table}[ht!]
\centering
\caption{Performance ratio}
\label{Tab:T3}
\begin{tabular}{|l|lll|lll|}
\toprule
 & \multicolumn{3}{c|}{Algorithm \ref{Alg:GA}} & \multicolumn{3}{c|}{Algorithm \ref{Alg:Fit}}\\
\textbf{n, $r_k$, T, R, mc} & \textbf{$C_{max}$}      & \textbf{WCT}       & \textbf{TWT}  & \textbf{$C_{max}$}      & \textbf{WCT}       & \textbf{TWT}     \\
\cmidrule{1-7}
(5,0,*,*,*)	&	1.02 (80)	&	1.01 (80)	&	1.02 (80)	&	1.02 (80)	&	1.01 (80)	&	1.07 (80)	\\
(5,1,*,*,*)	&	1.01 (80)	&	1.01 (80)	&	1.01 (80)	&	1.00 (80)	&	1.01 (80)	&	1.01 (80)	\\
(5,*,0.3,*,*)	&	1.01 (80)	&	1.01 (80)	&	1.01 (80)	&	1.01 (80)	&	1.01 (80)	&	1.07 (80)	\\
(5,*,0.6,*,*)	&	1.01 (80)	&	1.01 (80)	&	1.02 (80)	&	1.01 (80)	&	1.01 (80)	&	1.01 (80)	\\
(5,*,*,0.5,*)	&	1.01 (80)	&	1.01 (80)	&	1.01 (80)	&	1.01 (80)	&	1.01 (80)	&	1.06 (80)	\\
(5,*,*,2.5,*)	&	1.01 (80)	&	1.01 (80)	&	1.02 (80)	&	1.01 (80)	&	1.01 (80)	&	1.01 (80)	\\
(5,*,*,*,1)	&	1.01 (80)	&	1.01 (80)	&	1.01 (80)	&	1.01 (80)	&	1.01 (80)	&	1.00 (80)	\\
(5,*,*,*,2)	&	1.02 (80)	&	1.01 (80)	&	1.03 (80)	&	1.02 (80)	&	1.01 (80)	&	1.06 (80)	\\									
(15,0,*,*,*)	&	1.16 (20)	&	N/A	&	1.08 (45)	&	1.15 (20)	&	N/A	&	1.10 (45)	\\
(15,1,*,*,*)	&	1.02 (79)	&	1.02 (55)	&	1.02 (80)	&	1.01 (79)	&	1.01 (55)	&	1.10 (80)	\\
(15,*,0.3,*,*)	&	1.06 (50)	&	1.01 (24)	&	1.03 (76)	&	1.04 (50)	&	1.01 (24)	&	1.11 (76)	\\
(15,*,0.6,*,*)	&	1.04 (49)	&	1.02 (31)	&	1.05 (49)	&	1.03 (49)	&	1.01 (31)	&	1.08 (49)	\\
(15,*,*,0.5,*)	&	1.05 (51)	&	1.02 (28)	&	1.03 (59)	&	1.05 (51)	&	1.01 (28)	&	1.13 (59)	\\
(15,*,*,2.5,*)	&	1.04 (48)	&	1.02 (27)	&	1.04 (66)	&	1.03 (48)	&	1.01 (27)	&	1.08 (66)	\\
(15,*,*,*,1)	&	1.04 (55)	&	1.01 (35)	&	1.02 (64)	&	1.04 (55)	&	1.01 (35)	&	1.02 (64)	\\
(15,*,*,*,2)	&	1.05 (44)	&	1.03 (20)	&	1.06 (61)	&	1.04 (44)	&	1.02 (20)	&	1.18 (61)	\\										
(25,0,*,*,*)	&	N/A	&	N/A	&	1.13 (31)	&	N/A	&	N/A	&	1.28 (31)	\\
(25,1,*,*,*)	&	1.08 (78)	&	1.02 (8)	&	1.06 (77)	&	1.03 (78)	&	1.01 (8)	&	1.18 (77)	\\
(25,*,0.3,*,*)	&	1.08 (39)	&	1.02 (4)	&	1.05 (71)	&	1.02 (39)	&	1.01 (4)	&	1.20 (71)	\\
(25,*,0.6,*,*)	&	1.08 (39)	&	1.02 (4)	&	1.09 (37)	&	1.03 (39)	&	1.01 (4)	&	1.19 (37)	\\
(25,*,*,0.5,*)	&	1.09 (39)	&	1.01 (3)	&	1.07 (57)	&	1.03 (39)	&	1.01 (3)	&	1.17 (57)	\\
(25,*,*,2.5,*)	&	1.08 (39)	&	1.02 (5)	&	1.06 (51)	&	1.03 (39)	&	1.01 (5)	&	1.21 (51)	\\
(25,*,*,*,1)	&	1.04 (40)	&	1.02 (8)	&	1.05 (53)	&	1.01 (40)	&	1.01 (8)	&	1.08 (53)	\\
(25,*,*,*,2)	&	1.13 (38)	&	N/A	&	1.09 (55)	&	1.04 (38)	&	N/A	&	1.31 (55)	\\
\bottomrule
\end{tabular}
\end{table}

Table \ref{Tab:T3} presents a summary of the average performance ratios for every experimental factor of interest. This summary is segmented according to the number of jobs (n) and objective function under study. A ``*'' for an experimental factor denotes that all instances at all possible levels were combined. For example, the Table \ref{Tab:T3} row labeled (5,0,*,*,*) contains the average performance ratio of all instances with n = 5 jobs and no ready times (i.e., $r_k$ = 0), while the (25,*,*,2.5,*) rows contains the average performance ratio for all 25 job instances with due date range factor R = 2.5. For each objective function, the number given in parentheses denotes the number of optimal solutions that were found across the 80 instances analyzed. Finally, a Table \ref{Tab:T3} entry of ``N/A (0)'' denotes the case wherein no optimal solutions were found; therefore, no average performance ratio can be computed. 

Next performance measure is the heuristic ratio (HR) metric, which is defined as the ratio of the objective function value obtained by GA ($OF_{GA}$) for a problem instance to the non-optimal objective function value produced by the mathematical model in 7200 seconds ($OF_{7200}$) for the same problem instance. The average heuristic ratio for each experimental factor level by objective function is shown in Table \ref{Tab:T4}. A Table \ref{Tab:T3} entry of ``N/A (0)'' denotes the case wherein optimal solutions were found for all instances; therefore, no average heuristic ratio can be computed. The GA produced results that are 20\% above the time-limited MIP model solution for the makespan objective function. However, the 20\% above optimal performance is obtained in less than 5 minutes. One other observation that was obtained from Table \ref{Tab:T4} is that when the problem instance is small, the mathematical problem solved the instances to optimality for all scenarios and hence an HR is not available for those instances that had five jobs.


\begin{table}[ht!]
\centering
\caption{Heuristic ratio}
\label{Tab:T4}
\begin{tabular}{|l|lll|lll|}
\toprule
 & \multicolumn{3}{c|}{Algorithm \ref{Alg:GA}} & \multicolumn{3}{c|}{Algorithm \ref{Alg:Fit}}\\
\textbf{n, $r_k$, T, R, mc} & \textbf{$C_{max}$}      & \textbf{WCT}       & \textbf{TWT}  & \textbf{$C_{max}$}      & \textbf{WCT}       & \textbf{TWT}     \\
\cmidrule{1-7}
(15,0,*,*,*)	&	1.09 (60)	&	1.03 (80)	&	1.15 (34)	&	1.08 (60)	&	1.05 (80)	&	1.26 (34)	\\
(15,1,*,*,*)	&	1.01 (1)	&	1.04 (25)	&	N/A	&	1.01 (1)	&	1.05 (25)	&	N/A	\\
(15,*,0.3,*,*)	&	1.09 (30)	&	1.03 (56)	&	1.10 (3)	&	1.09 (30)	&	1.05 (56)	&	1.28 (3)	\\
(15,*,0.6,*,*)	&	1.08 (31)	&	1.03 (49)	&	1.15 (31)	&	1.08 (31)	&	1.06 (49)	&	1.25 (31)	\\
(15,*,*,0.5,*)	&	1.08 (29)	&	1.03 (52)	&	1.15 (20)	&	1.08 (29)	&	1.05 (52)	&	1.23 (20)	\\
(15,*,*,2.5,*)	&	1.09 (32)	&	1.03 (53)	&	1.14 (14)	&	1.09 (32)	&	1.06 (53)	&	1.29 (14)	\\
(15,*,*,*,1)	&	1.10 (25)	&	1.02 (45)	&	1.10 (15)	&	1.09 (25)	&	1.04 (45)	&	1.13 (15)	\\
(15,*,*,*,2)	&	1.08 (36)	&	1.04 (60)	&	1.18 (19)	&	1.08 (36)	&	1.06 (60)	&	1.35 (19)	\\										
(25,0,*,*,*)	&	1.16 (80)	&	1.02 (80)	&	1.16 (48)	&	1.17 (80)	&	1.10 (80)	&	1.36 (48)	\\
(25,1,*,*,*)	&	1.21 (2)	&	1.05 (72)	&	1.37 (3)	&	1.16 (2)	&	1.07 (72)	&	1.52 (3)	\\
(25,*,0.3,*,*)	&	1.15 (41)	&	1.04 (76)	&	1.19 (8)	&	1.17 (41)	&	1.08 (76)	&	1.27 (8)	\\
(25,*,0.6,*,*)	&	1.16 (41)	&	1.03 (76)	&	1.17 (43)	&	1.17 (41)	&	1.08 (76)	&	1.38 (43)	\\
(25,*,*,0.5,*)	&	1.16 (41)	&	1.03 (77)	&	1.14 (22)	&	1.17 (41)	&	1.07 (77)	&	1.31 (22)	\\
(25,*,*,2.5,*)	&	1.15 (41)	&	1.04 (75)	&	1.20 (29)	&	1.17 (41)	&	1.09 (75)	&	1.41 (29)	\\
(25,*,*,*,1)	&	1.19 (40)	&	1.03 (72)	&	1.14 (26)	&	1.19 (40)	&	1.07 (72)	&	1.27 (26)	\\
(25,*,*,*,2)	&	1.13 (42)	&	1.03 (80)	&	1.20 (25)	&	1.15 (42)	&	1.10 (80)	&	1.46 (25)	\\
\bottomrule
\end{tabular}
\end{table}

\subsection{Analysis}
{\color{blue}
From the results it is clear that the sets of instances, that have ready times, performed better than those sets of instances whose ready times are zero. This pattern is consistent irrespective of the number of jobs or the type objective function. The performance of the heuristic is also compared based on the two equipment scenarios mentioned in Table \ref{Tab:T2}. As the number of tools are reduced, the performance ratios is seen to be increased which is as expected. As the number of tools is reduced, the tightness of the resources increases which results in an increased value for average performance ratios. The parameters T and R have negligible impact on the average performance.

The results for the average heuristic ratio for set of instances with ready times and without ready times follow similar pattern as the average performance ratio. Similar to the performance ratio, heuristic ratio also seems to perform better for weighted completion time (WCT) and makespan ($C_{max}$) followed by the performance for total weighted tardiness (TWT). The GA fared better than the algorithm \ref{Alg:Fit} for the objective function WCT, TWT and most instances of $C_{max}$.
}

\section{Conclusion and Future Work}
A mixed-integer programming (MIP) formulation for the photolithography process with individual and cluster tool was developed for improved job scheduling. Due to this problem's complexity, a heuristic was developed to analyze our experimental cases. When comparing the the solution approaches, the MIP model provides better results but took a considerable amount of time. The heuristic approach achieved some good results in a very short span of time. The two heuristics' performance comparison showed the GA-based algorithm is efficient to schedule a photolithography process that involves larger number of jobs.  The heuristic method could be employed to scenarios when there is an unexpected machine downtime, shift changes, and changes in due dates. The heuristic seems to perform well for the objective function WCT and $C_{max}$. A better heuristic for TWT can be developed that can produce better results when compared to the current GA algorithm. Future work includes developing improved heuristic solutions to obtain better results by considering the job availability at every instant. The research can be extended to deal with finding a solution for minimizing other objectives like minimizing the total number of tardy jobs, minimizing maximum lateness, or to extend the research to investigate solutions for multiple objective problems.

\section*{References}

\bibliography{mybibfile}

\end{document}